\documentclass[twocolumn]{article}
\usepackage{amsfonts,eucal,bm,cite}

\begin{document}

\newtheorem{thrm}{Theorem}
\newtheorem{lem}{Lemma}
\newtheorem{ex}{Example}
\newtheorem{cor}[thrm]{Corollary}
\def\qed{\ \hfill\hspace*{\fill}
$\vbox{\hrule\hbox{\vrule height1.3ex\hskip1.3ex\vrule}\hrule}$
\hss\vskip\topsep\relax}

\def\bbR{\mathbb{R}}
\def\class{{\mathfrak A}}
\def\class{\ensuremath{\mathfrak A}}
\def\classP{\ensuremath{\mathfrak P}}
\def\classV{\ensuremath{\mathfrak V}}
\def\classR{\ensuremath{\mathfrak R}}
\def\classE{\ensuremath{\mathfrak E}}
\def\classF{\ensuremath{\mathfrak F}}

\def\norma{\ensuremath{\|\cdot\|}}
\def\span{\ensuremath{\mathrm{span}}}
\def\qcm{\ensuremath{\sigma}}
\def\ovm{\ensuremath{\chi}}
\def\absco{\ensuremath{\mathrm{absco}}}
\def\co{\ensuremath{\mathrm{co}}}
\def\set1N{\ensuremath{ \{1,2,\ldots ,N \}}}
\def\NbyN{\ensuremath{N\times N}}
\def\diag{\ensuremath{\mathrm{diag}}}
\def\length{\ensuremath{l}}

\def\proof{\textsc{Proof.\ }}

\title{\Large\textbf{Quasi-Controllability
and Estimates of Amplitudes of Transient Regimes in Discrete
Systems}\thanks{This work is supported in part by Russian
Foundation for Basic Research Grant 97-01-00692.}}

\author{\textsc{V.~Kozyakin, N.~Kuznetsov and
A.~Pokrovskii\thanks{Present address: Institute for Nonlinear
Science, Department of Physics, National University of Ireland,
University College, Cork, Ireland, e-mail:
alexei@peterhead.ucc.ie}}\\[5mm]
\textrm{Institute for Information Transmission Problems}\\
\textrm{Russian Academy of Sciences}\\
\textrm{19 Bolshoi Karetny lane, 101447 Moscow, Russia}
\\e-mail: kozyakin@iitp.ru\\[1cm]
{\normalsize\textbf{ABSTRACT}}\\[3mm]
\hskip -6pt\parbox{17.7cm}{\normalsize Families of
regimes for discrete control systems are studied possessing a
special quasi-controllability property that is similar to the
Kalman controllability property. A new approach is proposed to
estimate the amplitudes of transient regimes in
quasi-controllable systems. Its essence is in obtaining of
constructive a priori bounds for degree of overshooting in terms
of the quasi-controllability measure. The results are applicable
for analysis of transients, classical absolute stability problem
and, especially, for stability problem for desynchronized
systems.}}

\date{}

\maketitle

%

\section*{\centering\normalsize INTRODUCTION}
Currently, there are a growing number of cases in which systems
are described as operating permanently as if in the transient
mode. Examples are flexible manufacturing systems, adaptive
control systems with high level of external noises, so called
desynchronized or asynchronous discrete event systems
\cite{AKKK,KKKK,CHINA}. In connection with this, it is necessary
to ensure that the state vector amplitude satisfies reasonable
estimates within the whole time interval of the system
functioning including the interval of transient regime and an
infinite interval when the state vector is ``close to
equilibrium''. Emphasize,  that this necessity often contradicts
the usual desire to design a feedback which makes the system as
stable as possible. The reason is that the stability property
characterizes only the asymptotic behavior of a system and does
not take into account system behavior during the transient
interval. As a result, a stable system can have  large
overshooting or ``peaks'' in the transient process that can
result in complete failure of a system. First mentions about
systems with peak effects could be found in
\cite{BonYou68,BonYou70} and \cite{Mita76}. In
\cite{Kim81,Mita77,Mita78,Pol80b,Zei83} this effect was studied
for some classes of linear systems. As it was noted in
\cite{MiYo80,Izm87a,Izm88,Izm89,Kim81,Sha84,OlCi88} when the
regulator in feedback links is chosen to guarantee as large
degree of stability as possible then, simultaneously,
overshooting of the system state during the transient process
grows i.e., the peak effects  are getting more dangerous.

The above  papers were mainly concerned with continuous time
control systems because, in completely controllable and
observable discrete system it is possible to chose feedback
which turns to zero  the specter of the respective closed-loop
system. Nevertheless,  similar effects occurred when optimizing
asymptotic behaviour of  \emph{badly controllable or observable}
discrete time systems which arise in some applications.

Consider as the simplest, if trivial,  example the linear system
described by the relations
\begin{equation}\label{E0}
x_{n+1}=Ax_{n}+bx_{n}^{1},
\qquad
n = 0,1,2,\ldots .
\end{equation}
Here $x_{n}= \{x_{n}^{1},x_{n}^{2} \}\in\bbR^{2}$, vector
$b\in\bbR^{2}$ defines a feedback to be constructed and the
matrix
$$
A=\left(
\begin{array}{cc}
a & \varepsilon\\
\varepsilon & a
\end{array}
\right)
$$
depends on a small real parameter $\varepsilon$.

From the asymptotical point of view the best vector $b$ is
$(-2a,-(a^{2}+\varepsilon^{2})/\varepsilon)$ which makes
the eigenvalues of a closed system equal to zero. On the other
hand, for small $\varepsilon$ this vector $b$ is the most
dangerous at the first time step, because the closed system
(\ref{E0}) can be written for this $b$ as $x_{n+1}=A_{*}x_{n}$
where
$$
A_{*}=\left(
\begin{array}{cc}
-a & \varepsilon\\
a^{2}/\varepsilon & a
\end{array}
\right)
$$
has a big element $a^{2}/\varepsilon$ in the left bottom
corner.

There arises a general question if this kind of the peak effect
is connected only with poor controllability or observability of
the system? If an answer is positive, then the respective
quantitative estimates are of interest. Especially urgent such
estimates seems to be when a whole class of systems is examined
just as in problems of absolute stability or in desynchronized
systems. Another schemes of appearing peak  effects in discrete
systems see  in \cite{DIJ,IV}.

In this paper a new approach is developed presenting the means
to solve efficiently, for some classes of systems, the problem
of estimation the state vector amplitude within the whole time
interval. The key concept used is a quasi-controllability
property of a system that is similar to the Kalman
controllability property. The degree of quasi-controllability
can be characterized by a numeric value. The main result of the
paper is in proving the following:
\emph{if a quasi-controllable system is stable then the
amplitudes of all its state trajectories starting from the unit
ball are bounded by the value reciprocal of the
quasi-controllable measure.} Since the measure of
quasi-controllability can be easily computed, this fact becomes
an efficient tool for analysis of transients. It is shown also
that for quasi-controllable systems the properties of stability
or instability are robust with respect to small perturbation of
system's parameters. Some other results in this direction were
announced  in \cite{KP,KKP}. Note also that a concept similar to
quasi-controllability from algebraic point of view was
investigated in \cite{BW,DL}

\section*{\centering\normalsize
QUASI-CONTROLLABLE FAMILIES OF MATRICES}

The notion of \emph{quasi-controllability} of the system will be
introduced in this section. Degree of quasi-controllability will
be  estimated by some nonnegative value, the
\emph{quasi-controllability measure}. The basic  property of
quasi-controllability measure and some examples are also
discussed.

\subsection*{\normalsize Definitions and basic properties}
Let $\class =\{A_{1}$, $A_{2}$, \ldots , $A_{M}\}$ be a finite
family of real $N\times N$ matrices. The family $\class$ is said
to be \emph{quasi-controllable} one if no nonzero proper
subspace of $\bbR^{N}$  is invariant for all matrices from
$\class$. Evidently, a family $\class$ can be quasi-controllable
only if $M>1$.

Denote by ${\class}_{k}\ (k = 1,2,\ldots )$ the set of those
finite products of matrices from ${\class}\bigcup\{I\}$ which
contain no more that $k$ factors. Define ${\class}_{k}(x)$,
$x\in\bbR^{N}$ as the set of vectors $Lx$, with
$L\in{\class}_{k}$.  Denote by ${\co}(W)$ and ${\span}(W)$
respectively the convex and the linear hulls of the set
$W\subseteq\bbR^{N}$. Set also ${\absco}(W)={\co}(W\bigcup -W)$.
Let $\norma$ be a norm in $\bbR^{N}$; a ball in this norm of the
radius  $t$ centered at $0$ denote by $\mathbf{S}(t)$.

\begin{thrm}\label{L11}
Suppose that $p\ge N-1$. Then a  family $\class$ is
quasi-controllable if and only if ${\span}\{{\class}_{p}(x)\}
=\bbR^{N}$ for each nonzero $x\in\bbR^{N}$.
\end{thrm}

\proof
Let the family $\class$ be quasi-controllable and $x\in\bbR^{N}$
be a given nonzero vector. Introduce the sets ${\cal L}_{0} =
{\span}\{x\}$ and ${\cal L}_{k} = {\span}\{{\class}_{k}(x)\}$,
$k\ge 1$. Then
\begin{equation}\label{E11}
{\cal L}_{0}\subseteq{\cal L}_{1}\subseteq
\ldots \subseteq {\cal L}_{p}\subseteq\bbR^{N}.
\end{equation}
Therefore,
\begin{equation}\label{E12}
1\le\dim{\cal L}_{0}\le\dim{\cal L}_{1}\le\ldots\le\dim{\cal
L}_{p}\le N.
\end{equation}
On the other hand,
\begin{equation}\label{E13}
A_{i}{\cal L}_{j}\subseteq{\cal L}_{j+1}\qquad (A_{i}\in
{\class},\quad 0\le j\le p-1).
\end{equation}

If $\dim{\cal L}_{p}=N$ then ${\cal
L}_{p}={\span}\{{\class}_{p}(x)\} =\bbR^{N}$. If $\dim{\cal
L}_{p}<N$ then by (\ref{E12}) and the condition ${p\ge N-1}$,
the equality $\dim{\cal L}_{j}=\dim{\cal L}_{j+1}$ holds for
some ${j\in [0,p-1]}$. This and (\ref{E11}) imply ${\cal L}_{j}
= {\cal L}_{j+1}$. By the last equality and (\ref{E13}) the
subspace ${\cal L}_{j}$ should be invariant with respect to all
matrices from $\class$; due to quasi-controllability of the
family $\class$ this subspace coincides with $\bbR^{N}$. Hence,
${\cal L}_{j} = {\cal L}_{j+1} = \ldots = {\cal L}_{p} =
{\span}\{{\class}_{p}(x)\} =\bbR^{N}$.

Now suppose that  ${\span}\{{\class}_{p}(x)\} =\bbR^{N}$, but
the family ${\class}$ is not quasi-controllable. Then there
exists a nonzero proper subspace ${\cal L}\subset\bbR^{N}$ which
is invariant with respect to all matrices from $\class$. Then
the inclusion ${\span}\{{\class}_{p}(x)\}\subseteq{\cal L}$
holds for each ${x\in {\cal L}}$. Therefore,
${\span}\{{\class}_{p}(x)\}\neq\bbR^{N}$. This contradiction
proves the quasi-controllability of the family ${\class}$.  \qed

The value  {\qcm}$_{p}({\class})$ defined by
$$
{\qcm}_{p}({\class}) = \inf_{x\in\bbR^{N},\|x\|=1}
\sup\{t:\ \mathbf{S}(t)\subseteq{\absco}({\class}_{p}(x))\}
$$
is called \emph{$p$-measure of quasi-controllability} of the family
${\class}$ (with respect to the norm $\norma$).

\begin{thrm}\label{L12}
Let ${p\ge N-1}$. The family ${\class}$ is quasi-controllable if
and only if ${{\qcm}_{p}({\class})\neq 0}$.
\end{thrm}

\proof
Let ${{\qcm}_{p}({\class})\neq 0}$. Then the inclusion
$\mathbf{S}\{\|x\|{\qcm}_{p}({\class})\}\subseteq{\absco}\{{\class}_{p}
(x)\}$ holds for each nonzero ${x\in\bbR^{N}}$ and moreover
$\bbR^{N}={\span}\{{\class}_{p}(x)\}$. Therefore, by Theorem
\ref{L11} the family ${\class}$ is quasi-controllable.

Suppose now that the family $\class$ is quasi-controllable but
${\qcm}_{p}({\class})=0$. Then there exist $x_{n}\in\bbR^{N}$,
${\|x_{n}\|= 1}$ and $y_{n}\in{\absco}\{{\class}_{p}(x_{n})\}$
such that $y_{n}\rightarrow0$ and
$ty_{n}\not\in{\absco}\{{\class}_{p}(x_{n})\}$ for $t>1$.
Without loss of generality we can suppose that the sequences
$\{x_{n}\}$ and $\{y_{n}/\|y_{n}\|\}$ are convergent:
$x_{n}\rightarrow x$, $y_{n}/\|y_{n}\|\rightarrow z$.

By Theorem \ref{L11} the linear hull of the set
$\{{\class}_{p}(x)\}$ coincides with $\bbR^{N}$. Hence, there
exist matrices $L_{1}$, $L_{2}$, \ldots , $L_{N}\in{\class}_{p}$
such that the vectors $L_{1}x$, $L_{2}x$,\ldots , $L_{N}x$ are
linearly independent. Then the vectors $L_{1}x_{n}$,
$L_{2}x_{n}$, \ldots , $L_{N}x_{n}$ are also independent for all
sufficiently large $n$. It means that for any $n$ there exist
numbers $\theta^{(n)}_{1}$, $\theta^{(n)}_{2}$, \ldots,
$\theta^{(n)}_{N}$ such that
$\theta^{(n)}_{1}+\theta^{(n)}_{2}+\ldots+\theta^{(n)}_{N}=1$
and the vector
\begin{equation}\label{E14}
z_{n} =\sum^{N}_{i=1}\theta^{(n)}_{i}L_{i}x_{n}
\end{equation}
is collinear to  $y_{n}$ i.e., $z_{n}=\eta_{n}y_{n}$
$(\eta_{n}>0)$.

By definition, $z_{n}\in{\absco}\{L_{1}x_{n}$, $L_{2}x_{n}$,
\ldots , $L_{N}x_{n}\}\subseteq{\absco}\{{\class}_{p}(x_{n})\}$
and $ty_{n}\not\in{\class}_{p}(x_{n})$ for $t>1$. Therefore,
$\eta_{n}\le 1$. This and the condition $y_{n}\rightarrow 0$
imply $z_{n}\rightarrow 0$. Without loss of generality the
sequences $\{\theta^{(n)}_{1}\}$, $\{\theta^{(n)}_{2}\}$,
\ldots, $\{\theta^{(n)}_{N}\}$ can be supposed to be convergent
to some limits $\theta_{1}$, $\theta_{2}$, \ldots ,
$\theta_{N}$. Now, after transition to the limit in (\ref{E14}),
we get
$\theta_{1}L_{1}x+\theta_{2}L_{2}x+\ldots+\theta_{N}L_{N}x=0$
and $\theta_{1}+\theta_{2}+\ldots+\theta_{N}=1$. This
contradicts the linear independence of the vectors $L_{1}x$,
$L_{2}x$,
\ldots , $L_{N}x$, and the theorem is proved. \qed

The following theorem is useful when a family of matrices
depends on a parameter.
\begin{thrm}\label{L13}
Let ${p\ge N-1}$ and $N\times N$ matrices $A_{1}(\tau)$,
$A_{2}(\tau )$, \ldots, $A_{M}(\tau)$ be continuous at the point
$0$ with respect to the real parameter $\tau$. Suppose that the
family ${\class}(\tau) =\{A_{1}(\tau)$, $A_{2}(\tau)$,\ldots ,
$A_{M}(\tau)\}$ is quasi-controllable at the point $\tau=0$.
Then the system ${\cal A}(\tau )$ is quasi-controllable for all
sufficiently small $\tau$ and the function
${\qcm}_{p}({\class}(\tau))$ is continuous in $\tau$ at the
point  ${\tau=0}$.
\end{thrm}

\subsection*{\normalsize Examples}
Let $A$ be an $N\times N$ matrix and $b,c\in\bbR^{N}$. Consider
the family ${\class}={\class}(A,b,c)$ which consists of the
matrix $A$ and the matrix $Q=bc^{T}$ with entries
$q_{ij}=b_{i}c_{j}$, $i,j=1,\ldots,N$.
\begin{ex}\label{Ex10}
The family ${\class}(A,b,c)$ is quasi-controllable if and only
if the pair $(A,b)$ is completely controllable by Kalman and the
pair $(A,c)$ is completely observable.
\end{ex}

The following example is important in the theory of asynchronous
systems (see, e.g., \cite{AKKK,KKKK}). Consider an $N\times N$
scalar matrix $A = (a_{ij})$ and introduce the family
${\classP}_{1}(A) =\{A_{1}$, $A_{2}$,\ldots , $A_{N}\}$ by
equalities
$$
A_{i} =\left(\begin{array}{cccccc}
1 & 0 &\dots & 0 &\dots & 0\\
0 & 1 &\dots & 0 &\dots & 0\\
\vdots &\vdots &\ddots &\vdots &\ddots &\vdots\\
a_{i1}& a_{i2}&\dots & a_{ii}&\dots & a_{iN}\\
\vdots &\vdots &\ddots &\vdots &\ddots &\vdots\\
0 & 0 &\dots & 0 &\dots & 1\end{array}\right).
$$
Let the norm in $\bbR^{N}$ be defined as $\|x\| =|x_{1}| +
\ldots + |x_{N}|$.

\begin{ex}\label{Ex11}
The family ${\classP}_{1}(A)$ is quasi-controllable, if and only
if $1$ is not an eigenvalue of $A$ and the matrix $A$ is
irreducible. If ${\classP}_{1}(A)$ is quasi-controllable then
${\qcm}_{N}[{\classP}_{1}(A)]\ge {\alpha}{\beta}^{N-1}$ where
\begin{eqnarray*}
\alpha &=&{\frac{1}{2N}}\min\{\|(A-I)x\|:\|x\|=1\},\\
\beta &=&{\frac{1}{2}}\min\{|a_{ij}| : i\neq j, a_{ij}\neq0\}.
\end{eqnarray*}
\end{ex}

Consider an example influenced by discussions with
E.Kaszkurewich. Let the family ${\classV}(A)$ consists of the
matrices $A$, $D_{1}A$, $D_{2}A$,\ldots , $D_{N}A$. Here $A$ is
a scalar $N\times N$ matrix with the elements $a_{ij}$ and
$D_{i}$ $(i=1,2,\ldots ,N)$ are the diagonal matrices of the
form
$$
D_{i}= \diag \{d_{1i},d_{2i},\ldots ,d_{ii},\ldots ,d_{Ni}\},
$$
where $d_{ij}=1$ if $i\neq j$ and $d_{ij}=-1$ if $i=j$. Such
families are called the \emph{vertex families}. Let again $\|x\|
=|x_{1}| + \ldots + |x_{N}|$.

\begin{ex}\label{EX211}
The vertex family is quasi-controllable if and only if $0$ is
not an eigenvalue of $A$ and the matrix $A$ is irreducible. If
the family ${\classV}(A)$ is quasi-controllable then
${\qcm}_{N}[{\classV}(A)]\ge{\tilde\alpha}{\tilde\beta}^{N-1}$
where
\begin{eqnarray*}
\tilde\alpha& =&{\frac{1}{N}}\min\{\|Ax\|:\ \|x\|=1\},\\
\tilde\beta& =&\min\{|a_{ij}| :\  i\neq j, a_{ij}\neq 0\}.
\end{eqnarray*}
\end{ex}

\section*{\centering\normalsize QUASI-CONTROLLABILITY AND THE PEAK EFFECT}
This section contains the main results of the paper. There will
be investigated the influence of quasi-controllability on
stability, instability and transient processes of dynamical
systems which are generating by linear difference equations
\begin{equation}\label{E124}
x(n+1) = A(n)x(n)
\end{equation}
with varying coefficients. The conceptually simple and effective
method to estimate norms of solutions of difference equations
uniformly for all $n=0,1,2,\ldots$ on the whole interval of
existence, including the initial interval of transient mode as
well as the consequent infinite interval of ``asymptotic
behavior'', will be described.

\subsection*{\normalsize A priori estimate of overshooting measure}

Let  $\class$ be a family of $N\times N$ matrices. The
difference equation (\ref{E124}) is called
\emph{Lyapunov absolutely stable with respect to the family
$\class$}, if there exists $\mu<\infty$ such that for each
sequence $A(n)\in\class$ any solution $x(\cdot)$ of the
respective equation satisfies the estimate
\begin{equation}\label{E125}
\sup_{n\ge 0}\|x(n)\|\le\mu\|x(0)\|.
\end{equation}

The smallest $\mu$, for which estimates (\ref{E125}) hold is
called \emph{overshooting measure} of the equation (\ref{E124})
with respect to the family $\class$, and is denoted by
${\ovm}({\class})$.

\begin{thrm}\label{T11}
Let the equation {\rm (\ref{E124})} be Lyapunov absolutely
stable with respect to the  quasi-controllable family
${\class}$. Then the inequality
\begin{equation}\label{mainE}
{\ovm}({\class})\le\frac{1}{{\qcm}_{p}({\class})}
\end{equation}
holds for each $p\ge N-1$.
\end{thrm}

This assertion is the central result of the paper. The  proof is
relegated to the next subsection. Now we will discuss some ways
of its application. Clearly, the Lyapunov absolute stability of
the equation (\ref{E124}) is equivalent to the Lyapunov
stability of the difference inclusion
\begin{equation}\label{En124}
x(n+1) \in F_{\class}x(n).
\end{equation}
where the set-valued function $F_{\class}$ is defined by
$$
F_{\class}(x)=\co\{Ax:\ A\in\class\}.
$$
Inclusions of the form (\ref{En124}) embrace the usual systems
of the discrete absolute stability theory \cite{KrPo,Li,NC}. On
the other hand, the Lyapunov absolute stability follows from the
absolute stability of the respective system. Consequently, it is
possible to combine, when estimating  overshooting measure of
control systems, the classical methods of absolute stability
theory with Theorem \ref{T11}. An example of using this approach
will be presented later in discussion of applications to
desynchronized systems. Now let us give only a couple of the
simplest corollaries of Theorem~\ref{T11}.

\begin{cor}\label{2aC}
Let the family $\class$ be quasi-controllable and suppose that
the only uniformly bounded solution $\ldots,x_{-n},\ldots,
x_{-2},x_{-1},x_{0}$ of inclusion (\ref{En124}) is the zero
solution. Then for each $p\ge 1$ any solution  $x_{n}$,
$n=0,1,\ldots$ of the inclusion (\ref{En124}) satisfies the
inequality $\|x_{n}\|\le\|x_{0}\|{\qcm}_{p}^{-1}({\class})$.
\end{cor}

\proof
This corollary follows  from Theorem \ref{T11} and from the
principle of absence of bounded solutions in  absolute stability
problem (see, e.g., \cite{KrPo}). \qed

Consider now the difference equation
\begin{equation}\label{naE}
x_{n+1}=Ax_{n}+bu_{n},
\qquad
n=0,1,\ldots
\end{equation}
where $b\in\bbR^{N}$ and the scalars $u_{n}$ satisfy, for fixed
$\gamma>0$ and $c\in\bbR^{N}$, the inequality
$|u_{n}|\le\gamma|\langle c,x_{n}\rangle|$ (here
$\langle\cdot,\cdot\rangle$ denotes Euclidean inner product in
$\bbR^{N}$). Such equations are common in control theory (see,
e.g., \cite{Li}).

\begin{cor}\label{1aC}
Let the pair $(A,b)$ be completely controllable and the pair
$(A,c)$ be completely observable and suppose that
$\max_{|\omega|=1}\gamma|\langle c,(\omega
I-A)^{-1}b\rangle|<1$. Then for each $p\ge 1$ any solution
$x_{n}$, $n=0,1,\ldots$ of the equation (\ref{naE}) satisfies
the inequality
$\|x_{n}\|\le\|x_{0}\|{\qcm}_{p}^{-1}({\class}_{*})$ where
${\class}_{*}=\{A-\gamma\,bc^{T}, A+\gamma\,bc^{T}\}$.
\end{cor}

\proof
By virtue of the example \ref{Ex10} the class ${\class}_{*}$ is
quasi-controllable. So this corollary follows from Theorem
\ref{T11} and the circle criteria of absolute stability,
\cite{Li}.
\qed

\subsection*{\normalsize Proof of Theorem \protect\ref{T11}}
First of all, let us establish two auxiliary assertions. Denote by
${\classR}$ the set of all finite products of matrices from
$\class$. Define the \emph{length} $\length(R)$ of a matrix
$R\in\classR$ as a smallest quantity of factors
$A_{1},A_{2},\ldots ,A_{q}\in\class$ in the representation $R =
A_{1}A_{2}\ldots A_{q}$.

\begin{lem}\label{increase}
Let the family $\class$ be quasi-controllable and suppose that
the inequalities
\begin{equation}\label{EL1}
\|Rx_{*}\| > \mu\frac{1}{{\qcm}_{p}({\class})}\|x_{*}\|,\qquad \mu >1
\end{equation}
hold for some $x_{*}\in\bbR^{N}$ $(x_{*}\not=0)$, $p\ge N-1$,
$R\in\classR$. Then for any $x\in\bbR^{N}$, $x\not=0$, there
exists a matrix $R_{x}\in{\classR}$ such that
$\|R_{x}x\|\ge\mu\|x\|$, $\length(R_{x})\le\length(R)+p$.
\end{lem}

\proof
Fix an arbitrary $x\in\bbR^{N}$, $x\not=0$. The vector
${\qcm}_{p}({\class})x_{*}$ belongs to the set
$\absco\left\{{\class}_{p}
\left(\frac{\|x_{*}\|}{\|x\|}x\right)\right\}$
by the definition of the quasi-controllability measure.
Therefore, there exist scalars $\theta_{1}$, $\theta_{2}$,
\ldots , $\theta_{Q}$ with
\begin{equation}\label{E127}
\sum^{Q}_{i=1}|\theta_{i}|\le 1,
\end{equation}
and  matrices
$L_{1},L_{2},\ldots ,L_{Q}\in{\class}_{p}$
such that
$$
\sum^{Q}_{i=1}\theta_{i}\frac{\|x_{*}\|}{\|x\|}L_{i}x
={\qcm}_{p}({\class})x_{*}.
$$
Hence,
$$
\sum^{Q}_{i=1}\theta_{i}RL_{i}x
={\qcm}_{p}({\class})\frac{\|x\|}{\|x_{*}\|}Rx_{*},
$$
and, farther, by  (\ref{EL1}),
$$
\sum^{Q}_{i=1}\|\theta_{i}L_{i}Rx\|\ge\mu\|x\|.
$$
But then  (see (\ref{E127})) there exists $i$, $1\le i\le Q$,
such that the matrix $R_{x}=L_{i}R$ satisfies
$\|R_{x}x\|\ge\mu\|x\|$.

It remains to note that $\length(R_{x})\le\length(R)+p$, due to
inclusion $L_{i}\in{\class}_{p}$. So, the lemma is proved.
\qed

The equation (\ref{E124}) is said to be \emph{absolutely
exponentially unstable} (with the exponent $\lambda>1$) in the
family $\class$, if for some $\kappa>0$ and for each vector
$x\in\bbR^{N}$, $x\not=0$, there exists a sequence
$A(n)\in\class$, such that the solution $x(n)$ of the
corresponding equation (\ref{E124}) with the initial condition
$x(0)=x$ satisfies the estimate
\begin{equation}\label{Eexp}
\|x(n)\|\ge\kappa\lambda^{n}\|x(0)\|,\qquad n=0,1,2,\ldots .
\end{equation}

\begin{lem}\label{expo}
Let the family $\class$ be bounded and conditions of Lemma
\ref{increase} hold. Then the equation (\ref{E124}) is
absolutely exponentially unstable in the family $\class$.
\end{lem}

\proof
Fix an arbitrary vector $x\in\bbR^{N}$, $x\not=0$, and construct
an auxiliary sequence of vectors $\{z(m)\}$, $m=0,1,\ldots $, by
relations  $z(0)=x$ and
$$
z(m)=R_{z(m-1)}z(m-1),\qquad m=1,2,\ldots \ .
$$
Here $R_{z(m)}$ are the matrices from Lemma \ref{increase}. Then
by Lemma \ref{increase}
\begin{equation}\label{C1}
\|z(m)\|\ge\mu^{m}\|z(0)\|,\qquad m=0,1,2,\ldots \ .
\end{equation}

By definition, matrices $R_{z(m)}$, $m=0,1,\ldots $, can be
represented in the form
$$
R_{z(m)}=A_{m,l(m)},\ldots ,A_{m,2},A_{m,1},\qquad
A_{m,j}\in\class,
$$
where $l(m)$ is the length of $R_{z(m)}$. Denote by $\{A(n)\}$,
$n=0,1,\ldots $, the sequence of matrices
$$
A_{0,1},\ldots ,A_{0,l(0)},A_{1,1},\ldots ,A_{m,1},\ldots
,A_{m,l(m)},\ldots \ ,
$$
and consider the solution $x(n)$ of the respective equation
(\ref{E124}), with the initial condition $x(0)=x$. Then the
relations $x(q_{m})=z(m)$, $m=0,1,\ldots$, hold with $q_{0}=0$
and
$$
q_{m}=\sum^{m-1}_{i=0}l(i),\qquad m=1,2,\ldots .
$$
Estimates (\ref{C1}) imply
\begin{equation}\label{C2}
\|x(q_{m})\|\ge\mu^{m}\|x(0)\|,\qquad m=0,1,\ldots \ .
\end{equation}

Norms of matrices from $\class$ are uniformly bounded by the
conditions of the lemma and the estimates
\begin{equation}\label{C3}
q_{m}-q_{m-1}=l(m-1)\le K\qquad m=1,2,\ldots ,
\end{equation}
hold by Lemma \ref{increase}. Therefore, the inequality
(\ref{C2}), in a little bit weaker form, can be extended on the
positive integers $n$ from the interval $(q_{m-1},q_{m}]$:
\begin{equation}\label{C4}
\|x(n)\|\ge\nu\mu^{m}\|x(0)\|,
\end{equation}
where $\nu >0$, $q_{m-1}<n\le q_{m}$ and $m=0,1,\ldots$.

Now, taking into account that $q_{m}\le mK$ for  $m=0,1,\ldots$,
we obtain by virtue of  (\ref{C3}) that the inequalities
(\ref{C4}) for appropriate $\kappa>0$ ¨ $\lambda>1$ imply the
estimate (\ref{Eexp}). The lemma is proved. \qed

Let us return to and finish the proof of Theorem \ref{T11}.
Suppose that the theorem is false. Then there exists a sequence
of matrices $\{A(n)\in\class ,\ n=0,1,\ldots\}$ and a solution
$x(n)$ of the respective equation (\ref{E124}), such that
\begin{equation}\label{auE}
\|x(n_{0})\| > \frac{1}{{\qcm}_{p}({\class})}\|x(0)\|
\end{equation}
holds for some
$n_{0}\ge 1$, $p\ge{N-1}$.
The inequality (\ref{auE}) implies
$$
\|A(n_{0}-1)\ldots A(1)A(0)x(0)\| >
\frac{1}{{\qcm}_{p}({\class})}\|x(0)\|.
$$
Hence, by Lemma \ref{expo}, the equation (\ref{E124}) is
absolutely exponentially unstable in the family ${\class}$ and,
all the more, this equation is not Lyapunov absolutely stable
with respect to this family. Obtained contradiction proves the
theorem. \qed

\subsection*{\normalsize Application to desynchronized systems}\label{ex2S}

Recently much attention was paid to the development of methods
for the analysis of dynamics of multicomponent systems with
asynchronously interacting subsystems (see \cite{KKKK,CHINA} for
further references).  As examples we can mention the systems
with faults in data transmission channels, multiprocessor
computing and telecommunication systems, flexible manufacturing
systems and so on. It turned out that under weak and natural
assumptions systems of this kind possess such a strong property
as robustness. In applications the robustness is often treated
as reliability of a system with respect to perturbation of
various nature, e.g., drift of parameters, malfunctions or
noises in data transmission channels, etc.

Let us introduce basic notions of the desynchronized systems
theory. Consider a linear system $S$ consisting of $N$
subsystems $S_{1},S_{2}, \ldots ,S_{N}$ that interact at some
discrete instances $\{T^n\}$, $-\infty <n<\infty $.  The moments
of interaction may be chosen according to some deterministic or
stochastic law but generally they are not known in advance. Let
the state of each subsystem $S_{i}$ be determined within the
interval ${[T^{n},T^{n+1})}$ by a numerical value $x_{i}(n)$,
$-\infty <n<\infty $.

Let at each moment $T^{n}\in\{T^{k}:\ -\infty <k<\infty \}$ only
one of subsystems $S_{i}$ (say, with an index
$i=i(n)\in\{1,2,\ldots ,N\}$) may change its state and the law
of the state updating be linear:
$$
x_{i}(n+1)=\sum_{i=1}^{N}a_{ij}x_{j}(n),\qquad i=i(n).
$$
Consider the matrix $A=(a_{ij})$ and introduce for each
$i=1,2,\ldots ,N$ an auxiliary matrix $A_{i}$
($i$\emph{-mixture} of the matrix $A$) that is obtained from $A$
by replacing its rows with indexes $i\not =j$ with the
corresponding rows of the identity matrix $I$.  Then the
dynamics equation for the system $S$ can be written in the
following compact form:
\begin{equation}\label{unrel}
x(n+1)=A_{i(n)}x(n),\qquad -\infty <n<\infty .
\end{equation}
The system described above is referred to as \emph{the linear
desynchronized or asynchronous system}.
\begin{thrm}\label{Tdes}
Let $1$ be not an eigenvalue of $A$ and the matrix $A$ be
irreducible. Suppose that the desynchronized system  is Lyapunov
absolutely stable. Then
${\ovm}({\class})\le\left(\alpha{\beta}^{N-1}\right)^{-1}$.
\end{thrm}
The proof follows immediately from Theorem~\ref{T11} and
Example~\ref{Ex11}. Theorem~\ref{Tdes} is convenient to combine
with different methods of establishing stability of
desynchronized systems (see bibliography in \cite{AKKK}).

\section*{\centering\normalsize ROBUSTNESS OF INSTABILITY}

In this short conclusive section we will consider another
application of above methods to qualitative analysis of discrete
systems.

Consider the difference equation (\ref{E124}), where matrices
$A(n)$ belong to a family $\class(\tau)=\{A_{1}(\tau)$,
$A_{2}(\tau)$, \ldots, $A_{M}(\tau)\}$, which depends on a real
parameter $\tau$. Then a natural question arises about
dependence of properties of the equation (\ref{E124}) on the
parameter $\tau$.

\begin{thrm}\label{unst}
Let the family $\class(\tau)$ be quasi-controllable and
continuous at $\tau = 0$. Suppose that the equation (\ref{E124})
is not Lyapunov absolutely stable with respect to the family
$\class(0)$. Then the equation (\ref{E124}) is also not Lyapunov
absolutely stable (and, in fact, is absolutely exponentially
unstable) with respect to the family $\class(\tau)$ for all
sufficiently small $\tau$.
\end{thrm}

\proof
Let the equation (\ref{E124}) be not Lyapunov absolutely stable
with respect to the family $\class(0)$. Then there exist
matrices $A(n,\tau)\in\class(\tau)$, $n = 0,1,\ldots$, such that
at $\tau=0$ the solution of the respective  equation
(\ref{E124}) satisfies for some $n_{0}>0$ the inequality
$$
\|x(n_{0})\| > \frac{1}{{\qcm}_{N-1}[{\class}(0)]}\|x(0)\|.
$$
Therefore
$$
\|A(n_{0}-1,0)\ldots A(0,0)x(0)\| >
\frac{1}{{\qcm}_{N-1}[{\class}(0)]}\|x(0)\|.
$$
On the other hand, the matrices $\{A(n,\tau )\}$ are continuous
at the point $\tau =0$, so as the functions
${\qcm}_{N-1}[{\class}(\tau )]$ (see Theorem \ref{L13}). Hence
$$
\|A(n_{0}-1,\tau )\ldots A(0,\tau )x(0)\| >
\frac{1}{{\qcm}_{N-1}[{\class}(\tau)]}\|x(0)\|
$$
for all sufficiently small values of $\tau$. Then by Lemma
\ref{expo} the equation (\ref{E124}) is absolutely exponentially
unstable with respect to the class ${\class}(\tau)$.  The
theorem is proved.
\qed

In some situations the following corollary from Theorems
\ref{L13}, \ref{T11} and \ref{unst} is useful.

\begin{cor}\label{limstab}
Let a quasi-controllable family of matrices $\class$ =
$\{A_{1}$, $A_{2}$,\ldots , $A_{M}\}$ be the limit of families
$\class_{m}$ = $\{A_{1,m}$, $A_{2,m}$,\ldots , $A_{M,m}\}$.
Suppose that the equation (\ref{E124}) is Lyapunov absolutely
stable with respect to the families $\class_{m}$,
$m=1,2,\ldots$. Then this equation is Lyapunov absolutely stable
with respect to the family $\class$. More than that, the
families $\class_{m}$ are quasi-controllable and the  measures
of overshooting ${\ovm}(\class_{m})$ are uniformly bounded.
\end{cor}

The following two examples show that the previous corollary
turns to be false without the assumption about
quasi-controllability of the family ${\class}$.

\begin{ex}\label{limexp}
Consider the sequence of families $\classE_{m}$ = $\{E_{m}\}$,
each of which consists of the single matrix
$$
E_{m}=\left(\begin{array}{cc}
1-\frac{1}{m}&1\\0&1-\frac{1}{m}
\end{array}\right).
$$
Then the family $\classE=\lim_{m\to\infty}\{\classE_{m}\}$
consists of the matrix
$$
E=\left(\begin{array}{cc} 1&1\\0&1\end{array}\right),
$$
and is not quasi-controllable. Therefore the respective equation
(\ref{E124}) is not exponentially stable with respect to the
family $\classE$, notwithstanding that this equation is
exponentially stable with respect to the families $\classE_{m}$.
\end{ex}

\begin{ex}\label{ovmunbou}
Consider the sequence of the families $\classF_{m}=\{F_{m}\}$,
each of which consists of  the single matrix
$$
F_{m} =\left(\begin{array}{cc}
1-\frac{1}{m^2}&\frac{1}{m}\\
0&1-\frac{1}{m^2}
\end{array}\right).
$$
Then the respective limiting family $\classF$ includes only the
identity matrix $I$ and, therefore, is not quasi-controllable.
Evidently, the equation (\ref{E124}) is stable  with respect to
the family $\classF$, as well as with respect to the classes
$\classF_{m}$, $m=1,2,\ldots$. On the other hand, the measures
of overshooting ${\ovm}(\classF_{m})$ are not uniformly bounded.
\end{ex}

\end{document}